\documentclass[amsbook,12pt]{article}
\usepackage{epsfig, amsfonts,amssymb, amsmath, mathrsfs}

\usepackage[usenames]{color}
\usepackage{picinpar}

\usepackage{soul}

\usepackage{comment}

\newtheorem{theorem}{Theorem}[section]
\newtheorem{lemma}[theorem]{Lemma}

\newtheorem{proposition}[theorem]{Proposition}
\newtheorem{definition}[theorem]{Definition}
\newtheorem{corollary}[theorem]{Corollary}
\newtheorem{remark}[theorem]{Remark}

\title{Liouville's theorem and comparison results
for solutions of degenerate elliptic equations in exterior domains }

\author{Leonardo Bonorino \and Andr\'e Silva \and Paulo Zingano}

\date{}

\begin{document}

\maketitle

\begin{abstract}
A version of Liouville's theorem is proved for solutions
of some degenerate elliptic equations defined in $\mathbb{R}^n\backslash K$,
where $K$ is a compact set, provided the structure of this equation
and the dimension $n$ are related.
This result is a correction of a previous one established by Serrin,
since some additional hypotheses are necessary.
Theoretical and numerical examples are given.
Furthermore, a comparison result and the uniqueness of solution
are obtained for such equations in exterior domains.
\end{abstract}

\section{Introduction}

In this work we study some results established by J. Serrin in \cite{S}
concerning the classical Liouville's theorem and make some suitable corrections.
Furthermore,
we obtain a comparison principle and uniqueness result
for solutions in exterior domains.

According to Theorem 3 in \cite{S}, if
$u \in C^1(\mathbb{R}^n \backslash P)$
is a weak solution of
\begin{equation}
{\rm div} ( |\nabla v |^{p-2}A(|\nabla v|) \nabla v )  =  0
\quad {\rm in }\quad \mathbb{R}^n\backslash P,
\label{eqprincipal}
\end{equation}
bounded from below,
where $P = \{x_1, \dots, x_N\}$ and $A$ satisfies the following hypotheses,
\\ \null \quad (i) $p \ge n > 1$,
\\ \null \quad (ii) $A \in C([0,+\infty))$ with $A(0) > 0$,
\\ \null \quad (iii) $t^{\:\!p-1}A(t)$ is strictly increasing for $t > 0$,
\\ \null \quad (iv) $A(t)$ is bounded from above and also
bounded below for all $t \ge 0$ respectively
\\ \null \quad \quad \quad by positive constants $L >0$ and $\delta > 0$,
\\
then $u$ is constant.

This result (and its proof provided in \cite{S})
are correct only in the case $ p=n $;
simple counterexamples can be given if $p > n$.
Indeed,
if $A(t) \equiv 1$ and $ P=\{0\}$, then \eqref{eqprincipal}
is a $p$-\mbox{\small L}aplacian equation
that admits the radial solution
$$ u(x) = |x|^{\frac{p-n}{p-1}} \quad {\rm in } \quad \mathbb{R}^n \backslash \{0\},$$
which is bounded from below, but it is not constant.

This motivates the main goal of this paper,
which is to find additional hypotheses that
guarantee the Liouville's property and present some counterexamples when these assumptions fail.
In fact we establish some Liouville's extension in exterior domains
$\mathbb{R}^n \backslash K$, where $K$ is a compact set.
Related to this kind of question,
we also investigate the comparison principle in exterior domains
$\mathbb{R}^n \backslash K$
for a similar class of equations and
prove that if
$ u, v \in C^{0}(\overline{\mathbb{R}^n\backslash K}) \cap
  W^{1,p}_{\tt loc}(\mathbb{R}^n\backslash K) \cap
  L^{\infty}(\mathbb{R}^n\backslash K)$
are bounded weak solutions of
\begin{equation}
{\rm div} {\cal A} (x, \nabla u) = 0 \quad {\rm in} \quad  \mathbb{R}^n \backslash K
\label{eqSecundary2}
\end{equation}
satisfying $ u \le v$ on $\partial K$,
then $u \le v$ in $\overline{\mathbb{R}^n \backslash K}$.
From this, we get the uniqueness of bounded weak solutions
satisfying $u=f$ on $\partial K$,
where $f$ is a given boundary data.
To obtain the uniqueness and the comparison principle,
we need the following conditions on ${\cal A}$:
\\ \null \quad (v)
${\cal A} \in C(\overline{\mathbb{R}^n \backslash K}
\times \mathbb{R}^n;\mathbb{R}^n) \cap C^1(\overline{\mathbb{R}^n \backslash K}
\times \mathbb{R}^n \backslash \{0\}; \mathbb{R}^n)$
\\[10pt] \null \quad (vi)
${\cal A}(x,0) = 0 $ \quad for $x \in \mathbb{R}^{n}\backslash K$
\\[10pt] \null \quad (vii)
$\displaystyle \sum_{i,j\,=\,1}^n \left|
\frac{\partial {\cal A}_j(x,\eta)}{\partial \eta_i} \right|
\le \Gamma |\eta|^{p-2}$
\quad for $ x \in \mathbb{R}^n \backslash K $
and $ \eta \in \mathbb{R}^n \backslash \{0\} $
\\[10pt] \null \quad  (viii)
$\displaystyle \sum_{i,j\,=\,1}^n \frac{\partial {\cal A}_j(x,\eta)}{\partial \eta_i}
\xi_i \xi_j \ge \gamma |\xi|^2|\eta|^{p-2}$
\quad for $x \in \mathbb{R}^n \backslash K$,
$ \eta \in \mathbb{R}^n \backslash \{0\}$
and $ \xi \in \mathbb{R}^n,$
\\[5pt] where $\gamma$ and $\Gamma$ are suitable positive constants.
Although ${\cal A}$ can depend on $x$,
the problem studied in \eqref{eqprincipal} is not a particular case of
the one seen in \eqref{eqSecundary2},
since the conditions (v)-(viii) do not include some functions $A$
that satisfy (i)-(iv).

\

Several Liouville type theorems have been established for different equations.
For example, positive superharmonic functions in $\mathbb{R}^2$
must be constant (see \cite{PW})
and the same holds for positive solutions of
$\Delta_p u \le 0$ in $\mathbb{R}^n$ provided $p \ge n$ (see \cite{SZ} and \cite{MP}).
This result was extended by Mitidieri and Pohozaev to more general equations
in \cite{MP}, which was generalized by Filippucci in \cite{F}.
D'Alambrosio proved this result for a class of operators that includes
the sum of $p_i$-{\small L}aplacians with different degrees (see \cite{Da}).
Some interesting similar results were obtained for equations 
that have zero order terms as,
for instance, by Gidas and Spruck (\cite{GS}), by Bidaut-Veron (\cite{BV})
and by D'Alambrosio (\cite{Da}).
Other matter of interest is a Liouville type theorem for functions that are not defined
in all of $\mathbb{R}^n$, like the one investigated by Serrin in \cite{S}.
Another example is the work of Bidaut-Veron (\cite{BV})
that studies the case for exterior domains.

On the other hand, comparison results are also a subject of intense research.
For bounded domains, such results were obtained, for instance, by Guedda and Veron (\cite{GV}), Garc\'{\i}a-Meli\'an and Sabina de Lis (\cite{GL}), Cuesta and Tak\'a\v{c} (\cite{CT1}), Damascelli and Sciunzi (\cite{DS}), Roselli and Sciunzi (\cite{RS},\cite{S}), Galakhov (\cite{Ga2}) for equations involving the p-laplacian operator. For a more general class of quasilinear equations, some important comparison results were established by  Douglas et al (\cite{DDS}), Trudinger (\cite{Tr}), Tolksdorff (\cite{T}), Damascelli (\cite{D}), Pucci et al (\cite{PSZ}), Cuesta and Tak\`a\v{c} (\cite{CT2}), Lucia and Prashant (\cite{LP}), Pucci and Serrin (\cite{PS},\cite{PS2}), Alvino et al (\cite{ABM}), among others.

In the case that the domain is unbounded, 
comparison results are proved by Hwang (\cite{H}) for some 
quasilinear and mean curvature equations 
provided that some growth condition on the difference between the functions to be compared 
is satisfied. 
Uniformly elliptic linear equations are treated by 
Berestycki et al (\cite{BeHRo}) for general unbounded domains. 
They also investigate some semilinear equations in $\mathbb{R}^n$. 
For equations with the $p$-laplacian operator, 
comparison principles were established by Farina et al (\cite{FaMS}) and Sciunzi \cite{Sc} 
for ``narrow'' domains. 
For this kind of domain, some results are proved to the quasilinear case by Farina et al (\cite{FaMRS}). 
In \cite{Ga1}, Galakhov also obtained 
some comparison principle for quasilinear equations in domains that are the cartesian product between a bounded domain and some $\mathbb{R}^k$. 

Unlike the previous works, inspired in \cite{S}, 
this paper points out the importance of the condition $p \ge n$ 
to guarantee the validity of comparison principles for bounded solutions in exterior domains without requiring further assumptions on the domain. 
Also, no condition on the behavior of the solution at infinity is necessary. 
In other words, the boundary data determines the bounded solution 
in an exterior domain for $p \ge n$. 
(Such result is false in general for $p < n$.) 
From this we get some kind of Liouville's theorem in the sense that 
if the boundary data on an exterior domain is constant, 
then the constant function is the unique bounded solution for the problem.

In Section \ref{liouvillesection} we prove some Liouville's results
following the ideas of \cite{S}, making some necessary corrections.
The comparison principle and the uniqueness are proved for
bounded weak solutions in Section \ref{compUniquenessSection} for $ p \ge n $.
We show the existence of counterexamples in Section \ref{counterexample}
when some assumptions do not hold and
present an explicit counterexample computed numerically.

\section{Preliminary results}
\label{preliminary}

We need the following result,
which follows the same computation as in \cite{S}:

\begin{proposition}
There exists a family of radially symmetric weak solutions
$v_{a}(x)=v_{a}(|x|) \in  C^0(\mathbb{R}^n) \cap C^1(\mathbb{R}^n \backslash \{0\})$
of \eqref{eqprincipal} in $\mathbb{R}^n\backslash \{0\}$, for $p > n$,
such that
\\[5pt] {\rm (a)} $v_{a}(0)=0$ for any $a >0$;
\\ {\rm (b)} $r \to v_{a}(r)$ is strictly increasing and unbounded in $r$ for any fixed $a >0$;
\\ {\rm (c)} $\displaystyle \lim_{a \to 0} v_{a}(|x|) = 0$ for any $x \in \mathbb{R}^n$.
\label{radialbarrier}
\end{proposition}

{\sl Proof:}
To find increasing radially symmetric solutions of \eqref{eqprincipal},
we have to solve
\begin{equation}
\frac{d}{dr} \{(v')^{p-1} A(v')\} + \frac{n-1}{r} \{(v')^{p-1} A(v')\} = 0 ,
\label{radialDifferentialEquationFor1}
\end{equation}
which admits the solution
$$v_a(r)= \int_0^r \varphi^{-1}\left( \frac{a}{\tau^{n-1}} \right) d \tau ,$$
for any positive $a$, where $\varphi$, borrowed from \cite{S},
is defined by $\varphi(t)= t^{p-1} A(t)$.
Observe that from (ii) and (iii) we have that
$\varphi^{-1}$ is nonnegative, strictly increasing and continuous.
Moreover, condition (iv) implies that
\begin{equation}
 \left( \frac{a}{L}\right)^{1/(p-1)} \tau^{-\frac{n-1}{p-1}}
 \le \varphi^{-1}\left( \frac{a}{\tau^{n-1}}\right)
 \le \left( \frac{a}{\delta}\right)^{1/(p-1)} \tau^{-\frac{n-1}{p-1}}.
 \label{varphi-1inequality}
 \end{equation}
Then $v_a$ is well defined, $v_a(0)=0$ and
for any $r >0$ it holds
$$ 0 \le \lim_{a \to 0} v_a(r) \le \lim_{a \to 0}
\int_0^r \left( \frac{a}{\delta}\right)^{1/(p-1)} \tau^{-\frac{n-1}{p-1}} d \tau = 0,$$
proving (a) and (c).
Since $\varphi^{-1}$ is positive for $t >0$,
$v_a$ is strictly increasing.
The first inequality of \eqref{varphi-1inequality} implies that
$$ v_a(r) \ge \left( \frac{a}{L}\right)^{1/(p-1)} \left( \frac{p-1}{p-n} \right) r^{(p-n)/(p-1)},$$
showing that $v_a$ is unbounded,
concluding the result. \hfill $\square$

\begin{remark}
A similar result can be stated for $p=n$.
In this case, for a given $s >0$ there exists
a family of radially symmetric weak solutions
$\bar{v}_{a,s} = \bar{v}_a \in
C^1(\mathbb{R}^n \backslash B_s(0))$ in $\mathbb{R}^n \backslash B_s(0)$
that satisfies condition (b) and (c) of Proposition \ref{radialbarrier},
and satisfies also the following condition instead of (a):

(a') $\bar{v}_a(s) = 0$ for any $a$.
\label{auxiliaryPropositionForpn}
\end{remark}
For that, observe that
$$\bar{v}_a(r)= \int_s^r \varphi^{-1}\left( \frac{a}{\tau^{n-1}} \right) d \tau $$
is a solution of \eqref{radialDifferentialEquationFor1} with $p=n$, for $r \ge s$,
and satisfies condition (a').
As before, $\varphi^{-1}$ is nonnegative,
strictly increasing, continuous and
inequality \eqref{varphi-1inequality} is satisfied with $p=n$.
Hence, following the same steps, we can prove (b) and (c).

Consider the extension of $\bar{v}_a$ to $\mathbb{R}^n$ given by
$$ v_{a,s}(r) = v_a(r) = \left\{ \begin{array}{lc} \bar{v}_a(r) & {\rm if} \quad r \ge s\\[5pt]
0      & {\rm if} \quad r < s .\end{array} \right. $$
Observe that
$v_a \in C^0(\mathbb{R}^n) \cap C^1(\mathbb{R}^n \backslash \partial B_s(0))$
is a radially symmetric weak solution of \eqref{eqprincipal}
in $\mathbb{R}^n \backslash \partial B_s(0)$
that satisfies (a') and (c).
Condition (b) is also satisfied for $r \ge s$.

\

In order to develop the argument for comparison principle and uniqueness
we need Lemma 2.1 of \cite{D}, which in our case can be stated as follows:

\begin{lemma}
If ${\cal A}$ satisfies conditions (v)-(viii),
then there exist positive constants $c_1$ and $c_2$
depending on $p$, $\gamma$ and $\Gamma$ such that
\begin{equation}
|{\cal A}(x,\eta_2) - {\cal A}(x,\eta_1)| \le c_1 (|\eta_1|+|\eta_2|)^{p-2}|\eta_2 - \eta_1|
\label{damascelli1}
\end{equation}
and
\begin{equation}
[{\cal A}(x,\eta_2) - {\cal A}(x,\eta_1)]\cdot [\eta_2 -\eta_1] \ge c_2 (|\eta_1|+|\eta_2|)^{p-2}|\eta_2 - \eta_1|^2,
\label{damascelli2}
\end{equation}
for any $x \in \mathbb{R}^n \backslash K$ and $\eta_1,\eta_2 \in \mathbb{R}^n$
satisfying $|\eta_1|+|\eta_2| >0$,
where the dot stands for the scalar product in $\mathbb{R}^n$.
In particular,
\begin{equation}
|{\cal A}(x,\eta)| \le c_1 |\eta|^{p-1}
\label{damascelli3}
\end{equation}
\begin{equation}
{\cal A}(x,\eta) \cdot \eta \ge c_2 |\eta|^{p},
\label{damascelli4}
\end{equation}
for $x \in \mathbb{R}^n \backslash K$ and $\eta\in \mathbb{R}^n$.
Moreover, if $p \ge 2$, then
\begin{equation}
[{\cal A}(x,\eta_2) - {\cal A}(x,\eta_1)] \cdot
[\eta_2 -\eta_1] \ge c_2 |\eta_2 - \eta_1|^p,
\label{damascelli5}
\end{equation}
for $x \in \mathbb{R}^n \backslash K$ and $\eta_1,\eta_2 \in \mathbb{R}^n$.
\label{damascelli}
\end{lemma}
\

\section{ Liouville's result for solutions in exterior domains}
\label{liouvillesection}

To pursue the main result of this section
(Theorem 3.2),
we consider the following:

\begin{definition}
Given a compact set $K \subset \mathbb{R}^n$
and a function $u$ bounded in $ \mathbb{R}^n \backslash K$,
let $\,\overline{g}: \partial K \to \mathbb{R}$
and $\;\!\underline{g}: \partial K \to \mathbb{R}$
be defined by
$$ \overline{g} (x_0) =
\limsup_{\stackrel{\mbox{\scriptsize $x \to x_0$}}
{\mbox{$\scriptscriptstyle (\:\!\mbox{\scriptsize $x$} \,\in\, \mathbb{R}^{n}\backslash K)$}}}
\;\! u(x)
\quad {\rm and } \quad
\underline{g} (x_0) =
\liminf_{\stackrel{\mbox{\scriptsize $x \to x_0$}}
{\mbox{$\scriptscriptstyle (\:\!\mbox{\scriptsize $x$} \,\in\, \mathbb{R}^{n}\backslash K)$}}}
 u(x) $$
\end{definition}

\begin{theorem}
Let $K \subset \mathbb{R}^n$ be a compact set
and $u \in C^1(\mathbb{R}^n \backslash K)$
a weak solution of \eqref{eqprincipal} in $\mathbb{R}^n \backslash K$,
where conditions $(i)$-$(iv)$ are satisfied.
If $u$ is bounded and $p \ge n$, then
$$ I \le u(x) \le S, \qquad \forall \;\, x \in \mathbb{R}^{n} \backslash K $$
where $\,\displaystyle S = \sup_{x_0 \in \partial K} \overline{g}(x_0)$
and $\;\!\displaystyle I = \inf_{x_0 \in \partial K} \underline{g}(x_0)$.
\end{theorem}
{\sl Proof:}
Given $\varepsilon > 0$, we can prove that
$u < S + 3 \varepsilon$ in $\mathbb{R}^n \backslash K$.
For that, notice that from hypothesis we have that
for each $y \in \partial K$, there exists a ball $B_{\delta_y}(y)$
centered at $y$ of radius $\delta_y > 0$
such that
$$ u(x) \le \overline{g}(y) + \varepsilon \le S + \varepsilon
\quad {\rm in } \quad B_{\delta_y}(y) \cap \mathbb{R}^n \backslash K.$$
Since $\partial K$ is compact,
a finite number of these balls,
$B_{\delta_1}(y_1), \dots, B_{\delta_m}(y_m)$,
cover $\partial K$.
Defining
$V = (B_{\delta_1}(y_1) \cup  \dots \cup B_{\delta_m}(y_m)) \cap \mathbb{R}^n \backslash K$,
we have $u \le S +\varepsilon$ in $V$.

Let $\Omega \subsetneq V \cup K$ be a bounded open subset that contains $K$
such that $\partial \Omega \subset V$.
Hence, using that $V$ and $K$ are disjoint,
we have that $u$ is continuous on $\partial \Omega$ and
\begin{equation}
u \le S + \varepsilon \quad {\rm on} \quad \partial \Omega.
\label{psiupperestimate0}
\end{equation}
Let $x_0 \in V \backslash \overline{\Omega}$
and $R > 0$ such that $B_R(x_0) \supset \overline{\Omega}$.
Now
we show that $u < S + 3 \varepsilon$ in $B_R(x_0) \backslash K$,
defining $\psi_a : \mathbb{R}^n \to \mathbb{R}$ by
$$ \psi_a (x) = S + 2 \varepsilon + v_a(|x-x_0|), $$
where $v_a$ was introduced in Proposition \ref{radialbarrier} for $p >n$ and in Remark \ref{auxiliaryPropositionForpn} for $p=n$.

The idea is to use a comparison principle
to prove that $u \le \psi_a < S + 3 \varepsilon $
in $B_R(x_0) \backslash K$ for some suitable $a$.

First observe that $\psi_a$ is a weak solution of \eqref{eqprincipal}
in $\mathbb{R}^n\backslash \{x_0\}$ if $p >n$. For $p=n$, $\psi_a$ is a weak solution in
$\mathbb{R}^n \backslash B_s(0)$, where $s$ will be chosen later. In any case,
from (c) of Proposition \ref{radialbarrier} or Remark \ref{auxiliaryPropositionForpn},
$v_a(R) < \varepsilon$ if $a$ is chosen small enough.
Then, using that
$v_a(r)$ is nondecreasing,
we have
\begin{equation}
 \psi_a(x) \le S + 3 \varepsilon \quad {\rm for} \quad |x-x_0| < R.
 \label{psiupperestimate1}
\end{equation}
Furthermore, since $v_a$ is unbounded,
there exists $R_1 > R$ such that $ v_a(r) > M - S$ for $r \ge R_1$,
where $M=\sup  u(x)$.
Hence
\begin{equation}
 \psi_a(x) > M  \ge u(x) \quad {\rm for } \quad |x-x_0| \ge R_1.
\label{psilowerestimate1}
\end{equation}
Now,
using that $x_0 \in V$,
$u \le S+ \varepsilon $ in $V$ and $u$ is continuous,
it follows that
$u < S + 2\varepsilon$ in $B_{\rho}(x_0)$
for some small $\rho < R$.
(If $p=n$, we take $s=\rho$.)
Therefore
\begin{equation}
u < S + 2\varepsilon \le \psi_a \quad {\rm in } \quad B_{\rho}(x_0).
\label{psiupperestimate2}
\end{equation}
From this and \eqref{psilowerestimate1},
we have $u < \psi_a$ on $\partial B_{\rho}(x_0) \cup \partial B_{R_1}(x_0)$,
and from \eqref{psiupperestimate0}
we conclude that
$u \le S + \varepsilon < \psi_a$ on $\partial \Omega$.
Then
$$u <  \psi_a \quad {\rm on}  \quad  \partial \left( B_{R_1}(x_0)
\backslash (B_{\rho}(x_0) \cup \Omega )\right). $$
Hence, since $u$ and $\psi_a$ are solutions of \eqref{eqprincipal}
in $B_{R_1}(x_0) \backslash (B_{\rho}(x_0) \cup \Omega )$,
a comparison principle as in Theorem 2.4.1 of \cite{PS}
implies that
$u < \psi_a $ in $B_{R_1}(x_0) \backslash (B_{\rho}(x_0) \cup \Omega )$.
Then \eqref{psiupperestimate1}
implies that
$$ u < S + 3 \varepsilon \quad {\rm in} \quad B_{R}(x_0) \backslash (B_{\rho}(x_0) \cup \Omega ). $$
From \eqref{psiupperestimate2},
this inequality also holds in $B_{\rho}(x_0)$ and,
using that $u \le S+\varepsilon$ in $V$,
we have
$$  u < S + 3 \varepsilon \quad {\rm in} \quad B_{R}(x_0) \backslash K. $$
Since $R$ is arbitrary, we get $u < S + 3 \varepsilon$ in $\mathbb{R}^n  \backslash K$,
as stated in the beginning.
Hence $u \le S$. By an analogous argument, $u \ge I$,
completing the proof.
\hfill $\square$

\

Observe that $u$ is constant in the case $S=I$.
This is some version of the Liouville's theorem that can be stated as follows:
\begin{corollary}
Assuming the same hypotheses as in the previous theorem,
if $u$ can be extended continuously to the boundary of $K$
and its extension is constant on $\partial K$,
then $u$ is constant.
\label{uIsConstantIfSisEqualToI}
\end{corollary}

\begin{remark}
If $p < n$ the result is false. For instance, consider the
{\small L}aplacian equation in $\mathbb{R}^n$ for $n \ge 3$
and $K=\overline{B_{1}(0)}$.
This corresponds to equation \eqref{eqprincipal}
with $A(t)=1$ and $p=2$.
The function $u(x)= 1 - |x|^{2-n}$ is a bounded solution of the
{\small L}aplacian equation in $\mathbb{R}^n \backslash \overline{B_{1}(0)}$
such that $I=S=0$, but $u$ is not constant.
\end{remark}

From Corollary \ref{uIsConstantIfSisEqualToI}
we get a correct version of Theorem 3 in \cite{S}:

\begin{corollary}
Let $P=\{x_1, \dots, x_k\}$ be a finite set
and $u \in C(\mathbb{R}^n) \cap C^1(\mathbb{R}^n \backslash P) $
be a weak solution of \eqref{eqprincipal} in $\mathbb{R}^n \backslash P$,
where conditions $(i)$-$(iv)$ are satisfied.
Suppose that $u$ is bounded from below and from above,
$p \ge n$ and $u(x_1) = \dots = u(x_k)$.
Then $u$ is constant.
\end{corollary}

This result is not true if $p > n$ and
we do not suppose that $u$ is continuous in $P$,
$u(x_1) = u(x_2)= \dots = u(x_k)$, or $u$ is bounded from above.
For instance, if $u$ is not bounded from above,
Proposition \ref{radialbarrier} provides an example of a solution of \eqref{eqprincipal}
in $\mathbb{R}^n\backslash \{0\}$ that is continuous in all $\mathbb{R}^n$,
is bounded from below, but is not constant.
A counterexample for the case where $u$ is bounded
but condition $u(x_1) = u(x_2)= \dots = u(x_k)$ fails
is presented in the next section.

For $p = n$,
the continuity of $u$ on $\{x_1,\dots, x_k\}$,
$u(x_1)=\dots=u(x_k)$ and the boundedness from above
are unnecessary assumptions according to Theorem 3 of \cite{S}.

\

\section{Comparison principle and uniqueness for solutions in exterior domains}
\label{compUniquenessSection}

We can prove that the comparison principle for exterior domains holds
provided that the given boundary data and $A$
satisfy some conditions.
First we have to show some boundedness for $Du$.

\begin{theorem}
Let $ K $ be a compact set of $\mathbb{R}^n$
and $ u \in W^{1,p}_{\tt loc}(\mathbb{R}^n\backslash K)
\cap L^{\infty}(\mathbb{R}^n\backslash K) $
be a bounded weak solution of \eqref{eqSecundary2} in $\mathbb{R}^n\backslash K$.
Suppose that $ p \geq n $
and that ${\cal A}$ satisfies
$($v\/$)\:\!$-$\:\!($viii\/$)$.
If $ U \subset \mathbb{R}^n $ is an open set such that $ K \subset U $,
then $ |Du| \in L^p(\mathbb{R}^n\backslash U)$ and
$$ \| Du\|_{L^p( \mathbb{R}^n\backslash U)} \le C \sup_{\mathbb{R}^n\backslash K} |u|, $$
where $ C > 0 $ is a constant that depends on $p$, $U$, $K$
and the constants $c_1$ and $c_2$ of Lemma \ref{damascelli}.
\label{LpBoundnessOfDerivativeOfSolution}
\end{theorem}
{\sl Proof:}
We can assume without loss of generality that $U$ is bounded.
Let $R_0 >0$ such that $\overline{U} \subset B_{R_0}(0)$.
For $R > R_0$,
let $\psi=\psi_R \in C_0^{\infty}(\mathbb{R}^n)$
such that
\\$\bullet$ $0 \le \psi \le 1$
\\$\bullet$ $\psi(x) =1 $ for $x\in B_R(0) \backslash U$
\\$\bullet$ $\psi(x) = 0$ for $|x| \ge 2R$
\\$\bullet$ $\psi(x)=0 $ for $x \in V$,
where $V$ is a open set such that $K \subset V \subset \subset U$
\\$\bullet$ $|D\psi(x)| \le m$ for $x \in U\backslash V$,
where $m$ is a constant that depends on $U$ and $V$, but not on $R$
\\$\bullet$ $|D\psi(x)| \le 2/R$ for $ R < |x| < 2R$

\

Observe that the function $\varphi := \psi^p u$
belongs to $W^{1,p}_0(\mathbb{R}^n\backslash K)$.
Hence, using that
$u$ is a weak solution of \eqref{eqSecundary2}
and $\varphi$ an admissible test function,
we have
$$\int_E u \nabla \psi^p \cdot {\cal A}(x, \nabla u ) dx +
\int_E \psi^p \nabla u \cdot {\cal A}(x,\nabla u) dx = 0, $$
where $E \subset B_{2R}(0)\backslash V$
is the compact support of $\psi$.
Hence,
\eqref{damascelli3} and \eqref{damascelli4}
imply that
\begin{align*}
c_2 \int_E \psi^p |\nabla u|^p dx &\le
\int_E \psi^p \nabla u \cdot {\cal A}(x,\nabla u) dx \\[5pt]
& = - \int_E u \nabla \psi^p \cdot {\cal A}(x, \nabla u ) dx\\[5pt]
& \le \int_E p |u| \psi^{p-1}|\nabla \psi| | {\cal A}(x, \nabla u ) | dx \\[5pt]
 & \le c_1 p \; M \int_E \psi^{p-1}|\nabla \psi| |\nabla u|^{p-1} dx ,
\end{align*}
where $M=\sup_{\mathbb{R}^n \backslash K} |u|$.
Using H\"older inequality in the last integral,
we get
\begin{equation}
c_2 \int_E \psi^p |\nabla u|^p dx  \le  c_1 p \; M
\left(\int_E \psi^{p} |\nabla u|^{p} dx\right)^{\frac{p-1}{p}}
\left( \int_E |\nabla \psi|^p dx \right)^{\frac{1}{p}}.
\label{0thLpBoundOnU}
\end{equation}
Therefore,
\begin{equation}
 \int_E \psi^p |\nabla u|^p dx  \le
 \left(\frac{c_1 p \; M}{c_2} \right)^p  \int_E |\nabla \psi|^p dx .
\label{firstLpBoundOnU}
\end{equation}
From the hypotheses on $\psi$,
we have
\begin{align*}
 \int_E |\nabla \psi|^p dx &=
 \int_{U\backslash V} |\nabla \psi|^p dx +
 \int_{ \{ R \le |x| \le 2R\}}  |\nabla \psi|^p dx \\[5pt]
 & \le m^p \, |\,U\backslash V \;\!| \,+\,
 \left(\frac{2}{R}\right)^{\!p} \!\;\!
 |\,\{ R \le |x| \le 2R\;\!\} \;\!| \\[5pt]
 & = m^p \, |\,U\backslash V \;\!| \,+\,
 \left(\frac{2}{R}\right)^{\!p} \!
 \omega_n \,[\;\!(2R)^n - R^n]
\end{align*}
where $ |\;\!{\sf A}\;\!| $ stands for the
Lebesgue measure of ${\sf A}$
and $\omega_n$ is the volume of the unit ball.
Since $p \ge n$,
\begin{equation}
 \int_E |\nabla \psi|^p \, dx
 \,\le\,
 m^p \, |\, U\backslash V \;\!|
 \,+\,
 \frac{\omega_n \:\!2^{p+n}}{R^{p-n}}
 \,\le\,
 m^p \, |\,U\backslash V \;\!| \,+\,
\frac{\omega_n \:\!2^{p+n}}{R_0^{p-n}},
\label{LpEstimateForPsi}
\end{equation}
for $R > R_0$.
Observe that the expression in the right-hand side
does not depend on $R$.
It is a constant that we denote by $D$
and conclude from \eqref{firstLpBoundOnU}
that
$$ \int_E \psi^p \, |\nabla u|^p \,dx  \,\le\,
\left(\frac{c_1 \;\!p \;\! M}{c_2} \right)^{\!\!\;\!p} \! D .$$
Using that
$B_R(0) \backslash U \subset E$ and
$\psi_R=1$ in $B_R(0) \backslash U$,
we have
$$ \int_{B_R(0) \backslash U}  |\nabla u|^p \, dx
\,\le\,
\left(\frac{c_1 \;\!p \;\! M}{c_2} \right)^{\!\!\;\!p} \! D
\quad {\rm for \; any } \quad R > R_0. $$
Making $R \to \infty$, we conclude the result. \hfill $\square$

\

\begin{remark}
This result holds more generally. Assume that 
$ u \in W^{1,p}_{\tt loc}(\mathbb{R}^n\backslash K)
\cap L^{\infty}(\mathbb{R}^n\backslash K) $ is a weak solution of
$$ - {\rm div} {\cal A}(x,v) = f(x) \quad {\rm in} \quad \mathbb{R}^n \backslash K,$$
where $f\in L^1(\mathbb{R}^n \backslash U)$ for any $U$ such that $K \subset U$, 
${\cal A}$ satisfies (v)-(viii) and $p \ge n$. 
Then, we can prove that $|Du| \in L^p(\mathbb{R}^n\backslash U)$. 
Indeed using the same argument as in the previous result, 
we get the following inequality in place of \eqref{0thLpBoundOnU}:
$$ c_2 \!\int_E \psi^p \,|\nabla u|^p \,dx  \,\le\,  c_1 \;\!p \, M
\left(\int_E \psi^{p} \,|\nabla u|^{p} \,dx\right)^{\!\!\frac{p-1}{p}}
\!
\left( \int_E |\nabla \psi|^p \,dx \right)^{\!\!\;\!\frac{1}{p}} 
\!+\!\;\! \int_E \psi^p \,u \;\! f \,dx .$$
Hence, 
from \eqref{LpEstimateForPsi} and $|u| \le M$, 
we conclude that
$$ \int_E \psi^p \,|\nabla u|^p \,dx  \,\le\, D_1 \!
\left( \int_E \psi^p \, |\nabla u|^p \, dx \right)^{\!\!\frac{p-1}{p}} \!\!+\, D_2, $$
where $\displaystyle D_1 = \frac{c_1 \;\!p \, M }{c_2} 
\left( m^p \, |\,U\backslash V \;\!| \,+\,
\frac{\omega_n \:\!2^{p+n}}{R_0^{p-n}} \right)^{\!\!\;\!\frac{1}{p}} 
\!\!= \frac{c_1 \;\!p \, M }{c_2} D^{\frac{1}{p}}$ 
and $\displaystyle D_2 = \frac{M}{c_2} \!\int_E|f| \,dx$. \\
Therefore, 
we have at least one of the following inequalities:
$$ \int_E \psi^p \,|\nabla u|^p \,dx \,\le\, 
2\;\!D_1 
\left( \int_E \psi^p \,|\nabla u|^p \,dx \right)^{\!\!\frac{p-1}{p}} 
\quad 
{\rm or} 
\quad \int_E \psi^p \,|\nabla u|^p \,dx  \,\le\, 2\;\! D_2.$$
Analyzing both possibilities, 
we get
$$ \int_E \psi^p \,|\nabla u|^p \,dx  \,\le\, 
\max\,\{\;\!2^{\:\!p} D_1^{p},\, 2\;\! D_2\} 
\;\!\le\,
2^{\:\!p} D_1^p \,+\; 2\;\!D_2 
\;\!=\, 
\left(\frac{2\,c_1 \;\!p \, M }{c_2}\right)^{\!\!\;\!p} \! D 
\;+\;  
\frac{2\;\!M}{c_2} \! \int_E|f| \,dx. $$
Since $\psi_R=1$ in $B_R(0) \backslash U \subset E$ and 
$E \subset B_{2R}(0)\backslash V\!$, 
making $R \to \infty$
we obtain
$$ \left( \int_{\mathbb{R}^n \backslash U} |\nabla u|^p \,dx \right)^{\!\!\;\!1/p} 
\!\le\, \left(\frac{2\,c_1 \;\!p \, M }{c_2}\right) \!\;\!D^{1/p} 
\;\!+\,  
\left(\frac{2\;\!M}{c_2} \!
\int_{\mathbb{R}^n \backslash U} |f| \,dx\right)^{\!\!\;\!1/p} 
\!<\, \infty, $$
concluding the statement.
\label{ExtensionOfComparisonTheorem}
\end{remark}

Now we obtain a comparison result.

\begin{theorem}
Let $K$ be a compact set of
$\mathbb{R}^n$ and
$ u,v \in C^{0}(\overline{\mathbb{R}^n\backslash K}) \cap
  W^{1,p}_{\tt loc}(\mathbb{R}^n\backslash K) \cap
  L^{\infty}(\mathbb{R}^n\backslash K) $
be bounded weak solutions of \eqref{eqSecundary2}
in $\mathbb{R}^n\backslash K$.
Suppose that ${\cal A}$ satisfies
the conditions $($v\/$)\:\!$-$\:\!($viii\/$)$
and that $p \geq n$.
If $v \ge u$ on $\partial K$,
then
$$ v \ge u \quad {\rm in} \quad \mathbb{R}^n \backslash K.$$
\label{comparisonForSolutions}
\end{theorem}
\mbox{} \vspace{-1.250cm} \\
{\sl Proof\/}:
Let $\varepsilon > 0$.
Since $v-u + \varepsilon$ is continuous and
$v-u + \varepsilon \ge \varepsilon $ on $\partial K$,
there exists a bounded open set $U \supset K$
with smooth boundary
such that $ v-u + \varepsilon > 0$ in $U \backslash K$.
Let $R_0 \geq 1 $ and $\psi=\psi_R$, for $R > R_0$,
as described in the previous theorem.
Observe that
$$(v-u + \varepsilon)^- \;\!=\, \max \,\{-(v-u+\varepsilon),0\}
\:\!=\, 0
\quad {\rm in} \quad U \backslash K,$$
so that
$ \varphi : = \psi^p \;\!(v-u+\varepsilon)^-$ has a compact support
$ F = F_R \subset B_{2R}(0)\backslash U \subset \mathbb{R}^n\backslash K$.
Because
$ (v-u + \varepsilon)^- \in W^{1,p}_{\tt loc}(\mathbb{R}^n\backslash K)$
(see, for instance, Lemma 7.6 of \cite{GT}),
it follows that
$\varphi \in  W^{1,p}_0(\mathbb{R}^n\backslash K)$.
Let us also define
$ G_{\!\;\!\epsilon} \!:=
\{\;\! x \in \mathbb{R}^n\backslash K
\!\!\;\! : \;\!
v - u + \epsilon < 0 \,\}
\subset \mathbb{R}^n\backslash U \!\;\!$.
%
%

Since $u$ and $v$ are weak solutions
of \eqref{eqSecundary2},
taking the test function $\varphi$
we get
\begin{equation} \int_{\mbox{}_{\scriptstyle \mathbb{R}^{n} \backslash \;\!K}}
\hspace{-0.550cm}
\nabla \varphi \;\!\cdot [\;\! {\cal A}(x, \nabla v) -  {\cal A}(x, \nabla u)\;\!]
\, dx \,=\, 0.
\label{wealSolutionForUandV}
\end{equation}
Using that
(see e.g.\;\cite{GT})
$\nabla \varphi =\;\! \psi^p \;\!
\chi_{\mbox{}_{\scriptstyle G_{\epsilon}}} \!\!\nabla (u-v) \;\!+
p \;\!\;\! \psi^{p-1} (v-u+\varepsilon)^- \;\!\nabla \psi \;\!$
a.e.\;in $\mathbb{R}^{n} \backslash K $,
where
$ \:\!\chi_{\mbox{}_{\scriptstyle G_{\epsilon}}} \!\!\:\!$
is the characteristic function of
the set $ G_{\!\:\!\epsilon} $,
it follows that
\begin{align*}
\int_{\mbox{}_{\scriptstyle G_{\epsilon}}}
\hspace{-0.100cm}
\psi^p \, \nabla (v-u) &\cdot
[\;\! {\cal A}(x, \nabla v) -  {\cal A}(x, \nabla u)\;\!] \,dx \;\!\;\!= \\[5pt]
%
   &
   \int_{\mbox{}_{\scriptstyle \!\;\!B_{2R}(0) \backslash B_{R}(0)}}
   \hspace{-1.650cm}
   p \;\!\;\! \psi^{p-1} \;\!(v-u+\varepsilon)^- \,
   \nabla \psi
   \cdot [\;\! {\cal A}(x, \nabla v) -  {\cal A}(x, \nabla u) \;\!] \,dx.
\end{align*}
Hence,
from $p \ge n \ge 2$,
we can use \eqref{damascelli5}
to obtain
$$ c_2 \!
\int_{\mbox{}_{\scriptstyle G_{\epsilon}}}
\hspace{-0.100cm}
\psi^p \, | \nabla (u-v) |^p  \,dx \,
\le
\int_{\mbox{}_{\scriptstyle \!\;\!B_{2R}(0) \backslash B_{R}(0)}}
\hspace{-1.650cm}
p \;\!\;\! \psi^{p-1} \,(v-u+\varepsilon)^- \,
\nabla \psi \cdot
[\;\! {\cal A}(x, \nabla v) -  {\cal A}(x, \nabla u) \;\!] \,dx. $$
Then,
applying \eqref{damascelli3},
$$ c_2 \!
\int_{\mbox{}_{\scriptstyle G_{\epsilon}}}
\hspace{-0.100cm}
\psi^p \, | \nabla (u-v) |^p \;\!dx \,
\le \,
c_1 \;\! p \!\:\!
\int_{\mbox{}_{\scriptstyle \!\;\!B_{2R}(0) \backslash B_{R}(0)}}
\hspace{-1.650cm}
\psi^{p-1} \;\! M \, |\nabla \psi| \,
(\;\!|\nabla v |^{\:\!p-1} \!\;\!+\;\! |\nabla u|^{\:\!p-1})
\, dx, $$
where
$ M = \sup \,(v-u+\varepsilon)^- \le \;\!\sup|u| + \sup |v| + \varepsilon < \infty $.
Using the H\"older inequality,
this gives
\begin{align*}
 c_2 \!
 \int_{\mbox{}_{\scriptstyle G_{\epsilon}}}
 \hspace{-0.100cm}
 \psi^p  \, | \nabla (u-v) |^p \;\! dx \,
 & \le \,
 2 \, c_1 \;\!p\, M \!\!\;\!
 \int_{\mbox{}_{\scriptstyle \!\;\!B_{2R}(0) \backslash B_{R}(0)}}
 \hspace{-1.650cm}
 |\nabla \psi| \;\!\;\! \psi^{\:\!p-1} \,\!
 \max \,\{\:\!|\nabla u |, \;\!|\nabla v|\:\!\}^{p-1} \,dx \\[5pt]
 & \le\,
 2 \, c_1 \;\!p\, M \left(
 \int_{\mbox{}_{\scriptstyle \!\;\!B_{2R}(0) \backslash B_{R}(0)}}
 \hspace{-1.650cm}
 |\nabla \psi|^p \;\!dx\;\!\;\!
 \right)^{\!\!\!\;\!\frac{1}{p}}
 \!\:\!
 \left(
 \int_{\mbox{}_{\scriptstyle \!\;\!B_{2R}(0) \backslash B_{R}(0)}}
 \hspace{-1.650cm}
 \psi^{\:\!p} \;\!(\:\!|\nabla u |^{\:\!p} \!\;\!+\:\! |\nabla v|^{\:\!p}\:\!)
 \,dx
 \right)^{\!\!\!\!\;\!\frac{p-1}{p}} \!\!.
\end{align*}
Therefore,
since $ |\psi| \le 1 $ and $ p \geq 2 $,
\begin{align}
 \int_{\mbox{}_{\scriptstyle G_{\epsilon}}}
 \hspace{-0.100cm}
 \psi^{\:\!p} \;\! | \nabla (u-v) |^p \,dx
 \,
 &\le\;
 \frac{\:\!2 \;\!c_1 \:\!p\;\! M}{c_2} \;\!
 \left(
 \int_{\mbox{}_{\scriptstyle \!\;\!B_{2R}(0) \backslash B_{R}(0)}}
 \hspace{-1.650cm}
 |\nabla \psi|^{\:\!p} \,dx
 \,
 \right)^{\!\!\!\;\!\frac{1}{p}}
 \!\;\!
 \left(
 \int_{\mbox{}_{\scriptstyle \!\;\!B_{2R}(0) \backslash B_{R}(0)}}
 \hspace{-1.650cm}
 \psi^{\:\!p} \, (\:\!|\nabla u |^{\:\!p} \!\;\!+\:\! |\nabla v|^{\:\!p} \:\!)
 \, dx
 \right)^{\!\!\!\:\!\frac{p-1}{p}}
 \nonumber \\[5pt]
 &\le\;
 \frac{\:\!2 \;\! c_1 \:\!p\;\! M}{c_2} \;\!
 \left(
 \int_{\mbox{}_{\scriptstyle \!\;\!B_{2R}(0) \backslash B_{R}(0)}}
 \hspace{-1.650cm}
 |\nabla \psi|^{\:\!p} \, dx
 \,
 \right)^{\!\!\!\;\!\frac{1}{p}}
 \left(\;\!
 \|\;\!\nabla u\;\!
 \|_{\mbox{}_{\scriptstyle L^{p}(S_{R})}}
   ^{\mbox{}^{\scriptstyle p - 1}}
 \!+\,
 \|\;\!\nabla v \;\!
 \|_{\mbox{}_{\scriptstyle L^{p}(S_{R})}}
   ^{\mbox{}^{\scriptstyle p - 1}}
\:\!\right)
\!\:\!,
\label{comparisonTheoremIneq1}
\end{align}
where
$ S_{R} = B_{2R}(0) \backslash B_{R}(0) $.
Recalling that
$ p \geq n $
and
$ |\nabla \psi| \le 2/R$ in $ S_{R} $,
we have
\begin{equation}
 \int_{\mbox{}_{\scriptstyle \!\;\!B_{2R}(0) \backslash B_{R}(0)}}
 \hspace{-1.650cm}
 |\nabla \psi|^p \, dx \;
 \le\;
 \frac{2^p}{R^p}\:\omega_n \;\!R^n \:\![\;\!2^n \!-\!\;\! 1]
 \,\leq\,
 2^p \;\! [\;\!2^n \!-\!\;\! 1] \, \omega_n
 \label{comparisonTheoremIneq2}
\end{equation}
for all $ R > R_0$.
Therefore,
%
from \eqref{comparisonTheoremIneq1}
and
\eqref{comparisonTheoremIneq2},
this gives
$$
 \int_{\mbox{}_{\scriptstyle G_{\epsilon}}}
 \hspace{-0.100cm}
 \psi^p \;\!| \nabla (u-v) |^p \;\!dx
 \;\le\,
 \frac{\;\!8 \;\!c_1 \:\!p\;\! M}{c_2}
 \;\!\;\!
 \omega_{\mbox{}_{\scriptstyle n}}^{1/p} \,
 \bigl(\,
 \|\;\!\nabla u \;\!
 \|_{\mbox{}_{\scriptstyle L^{p}(S_{R})}}
   ^{\mbox{}^{\scriptstyle \:\!p - 1}}
 +\,
 \|\;\!\nabla v \;\!
 \|_{\mbox{}_{\scriptstyle L^{p}(S_{R})}}
   ^{\mbox{}^{\scriptstyle \:\!p - 1}}
 \bigr).
$$
By Theorem 4.1,
the right-hand side converges to $0$ as $R \to \infty$.
On the other hand,
 $$ 
 \int_{\mbox{}_{\scriptstyle G_{\epsilon}}}
 \hspace{-0.100cm}
 \psi^p_{R} \, | \nabla (u-v) |^p \, dx
 \;\to
 \int_{\mbox{}_{\scriptstyle G_{\epsilon}}}
 \hspace{-0.100cm}
 | \nabla (u-v) |^p \,dx
 \quad {\rm as}
 \quad R \to \infty, $$
so that,
for each $ \epsilon > 0 $,
we have
$$ 
 \int_{\mbox{}_{\scriptstyle G_{\epsilon}}}
 \hspace{-0.100cm}
 | \nabla (u-v) |^{\:\! p} \,dx \;\!\;\!=\;\!\;\! 0. $$
Since
$ \nabla ( v - u + \epsilon )^{-}
\!\;\!=\;\!
\chi_{\mbox{}_{\scriptstyle G_{\epsilon}}}
\!\!
\nabla (u - v) $
a.e.\;in $ \mathbb{R}^{n} \backslash K $,
it follows that
$$ 
 \int_{\mbox{}_{\scriptstyle \mathbb{R}^{n} \backslash K}}
 \hspace{-0.550cm}
 | \;\!\nabla (v - u + \epsilon)^{-} |^{\:\! p}\,dx
 \;\!\;\!=\;\!\;\! 0. $$
Thus,
$ (v - u + \epsilon)^{-} \!\:\!$
is constant in $ \mathbb{R}^{n} \backslash K $.
Recalling that
$ (v - u + \epsilon)^{-} \!\:\!=\:\!0 \:\!$
in $ U \backslash K $,
we therefore have
$ (v - u + \epsilon)^{-} \!\:\!=\:\! 0 $
everywhere in $ \mathbb{R}^{n} \backslash K $,
that is,
$$ u(x) \,\leq\, v(x) \;\!+\;\! \epsilon,
\qquad \forall \;\, x \in \mathbb{R}^{n} \backslash K. $$
Since
$ \epsilon > 0 $
is arbitrary,
this gives the result,
as claimed.
\hfill $\square$
%

\

A uniqueness result is a direct consequence of this theorem:

\begin{corollary}
Let $K$ be a compact set of $\mathbb{R}^n$
and $f$ be a continuous function on~$\partial K$.
Suppose that ${\cal A}$ satisfies the conditions
$($v\/$)\:\!$-$\:\!($viii\/$)$ for $p \geq n$.
If
$ u, v \in C(\overline{\mathbb{R}^n\backslash K}) \cap
  W^{1,p}_{\tt loc}(\mathbb{R}^n\backslash K)$
are bounded weak solutions of \eqref{eqSecundary2}
in $\mathbb{R}^n\backslash K$ satisfying $ u = f = v $ on $\partial K$,
then
$ u = v $ in $\mathbb{R}^n \backslash K$.
\label{uniquenessOfSolutions}
\end{corollary}

\begin{corollary}
Let $K$ be a compact set of $\mathbb{R}^n$ and
$f$ be a continuous function on $\partial K$.
Suppose that ${\cal A}$ satisfies the conditions
$($v\/$)\:\!$-$\:\!($viii\/$)$ for $ p \geq n $.
If
$ u \in C(\overline{\mathbb{R}^n\backslash K}) \cap
  W^{1,p}_{\tt loc}(\mathbb{R}^n\backslash K)$
is a bounded weak solution of \eqref{eqSecundary2}
in $\mathbb{R}^n\backslash K$
satisfying $ u = f $ on $\partial K$,
then
$$ \min f \,\le\,  u  \,\le\, \max f
\quad {\rm in } \quad \mathbb{R}^n\backslash K. $$
\label{limitationOfuSecondResult}
\end{corollary}
{\sl Proof:}
Observe that $v\equiv \max f$ is a solution of \eqref{eqSecundary2}
and satisfies $v \ge u$ on $\partial K$.
Then, from Theorem \ref{comparisonForSolutions},
we get $u \le \max f$.
Similarly, $u \ge \min f$. \hfill $\square$

\begin{remark}
Observe that if $u,v \in
  W^{1,p}_{\tt loc}(\mathbb{R}^n\backslash K)$ satisfy
$$ -{\rm div}{\cal A}(x,u) \le -{\rm div}{\cal A}(x,v) \quad {\rm in } \quad \mathbb{R}^n\backslash K$$
in the weak sense, then we have
$$ \int_{\mbox{}_{\scriptstyle \mathbb{R}^{n} \backslash \;\!K}}
\hspace{-0.550cm}
\nabla \varphi \;\!\cdot [\;\! {\cal A}(x, \nabla v) -  {\cal A}(x, \nabla u)\;\!]
\, dx \,\ge\, 0, $$
instead of \eqref{wealSolutionForUandV}, for any nonnegative $\varphi \in W^{1,p}_0(\mathbb{R}^n\backslash K)$. Hence, if we assume that $|D u|, |Dv| \in L^p(\mathbb{R}^n\backslash U)$ for any open set $U$ such that $K \subset U$ and $u,v \in C^{0}(\overline{\mathbb{R}^n\backslash K}) \cap
  W^{1,p}_{\tt loc}(\mathbb{R}^n\backslash K) \cap
  L^{\infty}(\mathbb{R}^n\backslash K)$, following the same steps as in Theorem \ref{comparisonForSolutions}, we obtain a similar comparison principle:
\end{remark}

\begin{theorem}
\label{comparisonTheoremForNonhomegeneousEquation}
Let $K$ be a compact set of
$\mathbb{R}^n$ and
$ u,v \in C^{0}(\overline{\mathbb{R}^n\backslash K}) \cap
  W^{1,p}_{\tt loc}(\mathbb{R}^n\backslash K) \cap
  L^{\infty}(\mathbb{R}^n\backslash K) $
be bounded functions such that
$$ -{\rm div}{\cal A}(x,u) \le -{\rm div}{\cal A}(x,v) \quad {\rm in } \quad \mathbb{R}^n\backslash K$$
in the weak sense. 
Suppose that $|D u|, |Dv| \in L^p(\mathbb{R}^n\backslash U)$ 
for any open set $U$ such that $K \subset U$,
that ${\cal A}$ satisfies
the conditions $($v\/$)\:\!$-$\:\!($viii\/$)$
and that $p \geq n$.
If $v \ge u$ on $\partial K$,
then
$$ v \ge u \quad {\rm in} \quad \mathbb{R}^n \backslash K.$$
\end{theorem}

\begin{corollary}
Let $K$ be a compact set of
$\mathbb{R}^n$ and
$ u,v \in C^{0}(\overline{\mathbb{R}^n\backslash K}) \cap
  W^{1,p}_{\tt loc}(\mathbb{R}^n\backslash K) \cap
  L^{\infty}(\mathbb{R}^n\backslash K) $
be bounded functions such that
$$ -{\rm div}{\cal A}(x,u) = f_1(x) \le f_2(x) = -{\rm div}{\cal A}(x,v) \quad {\rm in } \quad \mathbb{R}^n\backslash K$$
in the weak sense, where $f_1,f_2 \in  L^1(\mathbb{R}^n\backslash U)$ for any open set $U$ such that $K \subset U$.
Suppose that ${\cal A}$ satisfies
the conditions $($v\/$)\:\!$-$\:\!($viii\/$)$
and that $p \geq n$.
If $v \ge u$ on $\partial K$,
then
$$ v \ge u \quad {\rm in} \quad \mathbb{R}^n \backslash K.$$
\end{corollary}
{\sl Proof:} Since $f_1$ and $f_2$ satisfy the hypothesis of Remark \ref{ExtensionOfComparisonTheorem}, it follows that $|Du|$ and $|Dv|$ are in $L^p(\mathbb{R}^n \backslash U)$ for any open set $U$ such that $K \subset U$. Therefore, from Theorem \ref{comparisonTheoremForNonhomegeneousEquation}, we get the result. \hfill $\square$
\begin{remark}
Observe that if $A$ satisfies conditions (i), (ii), (iv) and
\\ \null \quad (iii)'  $A(t)$  is nondecreasing in $[0, +\infty)$,
\\ then ${\cal A}(x,v) = |v|^{p-2}v A(|v|)$ satisfies conditions (v)-(viii). Hence all results in this section holds for $|v|^{p-2}v A(|v|)$ if conditions (i),(ii),(iii)' and (iv) are satisfied.
\end{remark}

\section{Examples of non constant bounded solutions}
\label{counterexample}

In this section
we build an example of weak solution of \eqref{eqprincipal}
in $\mathbb{R}^n \backslash \{x_1, \dots,  x_k\}$
that is bounded, belongs to $C(\mathbb{R}^n) \cap C^1(\mathbb{R}^n \backslash \{x_1, \dots, x_k\})$, but is not constant.

Indeed, for
$p > n$, $P = \{x_1, \dots, x_k \} \subset \mathbb{R}^n$
and $m_1, \dots, m_k \in \mathbb{R}$,
we show that the problem
\begin{equation}
\left\{ \begin{array}{rcl} \Delta_p u & = & 0 \quad {\rm in } \quad \mathbb{R}^n \backslash{P}\\[5pt]
                       u(x_i) & = & m_i \quad {\rm for } \quad i\in\{1, \dots, k\} \\[5pt]
\end{array} \right.
\label{plaplacianequationonRnminusP}
\end{equation}
has a bounded weak solution
in $C(\mathbb{R}^n) \cap C^1(\mathbb{R}^n \backslash P)$.

\begin{proposition}
Problem \eqref{plaplacianequationonRnminusP} has a weak solution
in $C(\mathbb{R}^n) \cap C^1(\mathbb{R}^n \backslash P)
\cap L^{\infty}(\mathbb{R}^{n})$.
\end{proposition}
{\sl Proof:}
First,
let $m= \min \;\!\{ m_1, \dots , m_k\}$,
$M= \max \;\!\{ m_1, \dots , m_k\}$ and $r_0 > 0$
such that the balls $B_{r_0}(x_1), \dots , B_{r_0}(x_k)$
are disjoint.
Consider also $R_0 > 0$
such that $\overline{B}_{r_0}(x_i) \subset B_{R_0}(0)$
for all $i \in \{1, \dots, k\}$,
and define
$$\Omega_{r,R} = B_{R} (0) \backslash  \overline{B_{r}(x_1) \cup \dots \cup B_{r}(x_k)}  \quad {\rm for } \quad r < r_0 \quad {\rm and} \quad R > R_0.$$
Using standard techniques,
we can show that there exists a weak solution $u_{r,R}$
of
\begin{equation}
\left\{ \begin{array}{rcl} \Delta_p u & = & 0 \quad {\rm in } \quad \Omega_{r,R}\\[5pt]
                       u & = & m_i \quad {\rm on } \quad \partial B_{r}(x_i) \\[5pt]
                       u & = & \frac{m+M}{2} \quad {\rm on } \quad \partial B_{R}(0)
\end{array} \right.
\end{equation}
We also know that $u_{r,R} \in C^{1,\beta}(\overline{\Omega}_{r,R})$
for some $\beta \in (0,1)$
(for instance, see \cite{E} or Theorem 6.20 of \cite{MZ}).
By the maximum principle, $m \le u_{r,R} \le M$ in $\Omega_{r,R}$.
For each $i \in \{1, \dots, k\}$,
define
$$ \psi_i^+(x) = m_i + v_{a_i}(|x-x_i|) \quad {\rm and } \quad \psi_i^-(x) = m_i-v_{a_i}(|x-x_i|), $$
where $v_{a_i}$ is as introduced in Proposition \ref{radialbarrier},
associated to the $p$-{\small L}aplacian problem,
and $a_i$ is chosen such that
$$ \psi_i^+(x) \ge M \quad {\rm if} \quad |x - x_i| \ge r_0
\quad {\rm and} \quad
\psi_i^-(x) \le m \quad {\rm if} \quad |x - x_i| \ge r_0.$$
Observe that $\psi_i^+$ and $\psi_i^-$ are weak solutions of
$\Delta_p v=0$ in $\mathbb{R}^n \backslash \{x_i\}$,
radially symmetric with respect to $x_i$
and $\psi_i^+(x_i)=\psi_i^-(x_i) = m_i$.

Since $\psi_i^+$ and $u$ are weak solutions of $\Delta_p v=0$ in $\Omega_{r,R}$,
$u_{r,R} \le M \le \psi_i^+$ on $\partial B_{r_0}(x_i)$
and $u_{r,R} = m_i < \psi_i^+$ on $\partial B_{r}(x_i)$,
the comparison principle implies that
$u_{r,R} \le \psi_i^+$ in the annulus $B_{r_0}(x_i) \backslash B_r(x_i)$.
Hence, using that $u_{r,R}(x) \le M < \psi_i^+(x)$ if $|x-x_i| > r_0$
(recall that $v_{a_i}(|x-x_i|)$ increases as $|x-x_i|$ increases and $\psi_i^+(r_0) \ge M$),
we have
$$ u_{r,R} \le \psi_i^+ \quad {\rm in } \quad \Omega_{r,R}
\quad {\rm for \; any} \quad i \in \{1, \dots, k\}.$$
In the same way,
$$ u_{r,R} \ge \psi_i^- \quad {\rm in } \quad \Omega_{r,R}
\quad {\rm for \; any} \quad i \in \{1, \dots, k\}.$$
For $R$ fixed,
let $(r_j)$ be a decreasing sequence converging to $0$
and $w_j := u_{r_j,R}$.
Given a compact set $K$ of $B_R(0) \backslash P$,
we have $K \subset \Omega_{r_j,R}$ for $j$ large.
Moreover, since $m \le w_j \le M$,
Theorem 1.1 of \cite{LU} (page 251) implies that
$$ |w_j(y)-w_j(x)| \le C |y-x|^{\beta} \quad {\rm for \; any} \quad x,y \in K,$$
for some $\beta=\beta(n,p) \in (0,1)$ and $C>0$,
where $C$ depends only on $n$, $p$, $m$, $M$ and $dist(K,\partial (B_R(0) \backslash P))$.
Hence, Arzel\`a-Ascoli's Theorem guarantees that
some subsequence of $w_j$, that we rename by $w_k$,
converges uniformly to some continuous function in $K$.
Considering a increasing sequence of compact subsets $K_{\ell}$
such that $B_R(0) \backslash P = \bigcup K_{\ell}$
and using a diagonal process,
we can find a continuous function $w:B_R(0) \backslash P \to \mathbb{R}$
and a subsequence $w_k$ that converges uniformly to $w$
on compact sets of $B_R(0) \backslash P$.

We define $u_R:=w$.
Observe that the uniform convergence of $w_k$ to $u_R$ in compacts
implies that $u_R$ is a weak solution of $\Delta_p v=0$ in $B_R(0) \backslash P$.
Hence, $u_R \in C^{1,\beta'}_{\tt loc}( B_R(0) \backslash P)$
for some $\beta' \in (0,1)$.
Then $u_R \in C^1(B_{R}(0)\backslash P)$ and,
using that
$$\psi_i^- \le w_k \le \psi_i^+
\quad {\rm for \; any } \quad  1\le i \le n \quad {\rm and} \quad k \in \mathbb{N},$$
we have that $\psi_i^- \le u_R \le \psi_i^+$
for any $i \in \{1, \dots ,n\}$.
Therefore,
$\lim_{x \to x_i} u_R(x) =m_i$.
Thus $u_R$ can be extended continuously to $B_R(0)$,
satisfying
\begin{equation}
 u_R(x_i) = m_i \quad {\rm for } \quad i \in \{1, \dots, n\}.
 \label{valueOfuRinxi}
\end{equation}
Finally,
taking a sequence $R_k \to \infty$, using a compactness argument
and a diagonal process,
we get, as before, a subsequence $u_{R_k}$
converging to some function
$u \in C(\mathbb{R}^n) \cap C^1(\mathbb{R}^n \backslash P)$
that is a weak solution of $\Delta_p v =0$ in $\mathbb{R}^n \backslash P$.
From \eqref{valueOfuRinxi}, we get $u(x_i)=m_i$,
proving the existence of solution to \eqref{plaplacianequationonRnminusP}.
\hfill $\square$ \\
\begin{remark}
This proposition also holds if we replace the $p$-{\small L}aplacian equation
by equation \eqref{eqprincipal},
provided that, besides the conditions (i)-(iv), 
$A$ also satisfies (v)-(viii)
to guarantee the $C^{1,\alpha}$ regularity of weak solutions.
\end{remark}

In the 2-D case,
we illustrate in Fig.\;1
the solution $ u = u(x,y) $
to Problem (7)
as constructed in the proof of Proposition 4.1
in the particular case $ n = 2 $, $ p = 4 $,
$ P = \{\;\!(-1,0), \;\!(1,0)\;\!\} $,
with $ \;\!m_1 = -\;\!1 $,
$ m_2 = 1 $.
Note the expected symmetries
($u(x,y) $ is odd in $x$,
and even in $y$)
and the asymptotic zero value
in the far field.

\mbox{} \vspace{-0.350cm} \\
\mbox{} \hspace{+0.000cm}
\begin{minipage}[t]{14.500cm}
\includegraphics[keepaspectratio=true, scale=1]{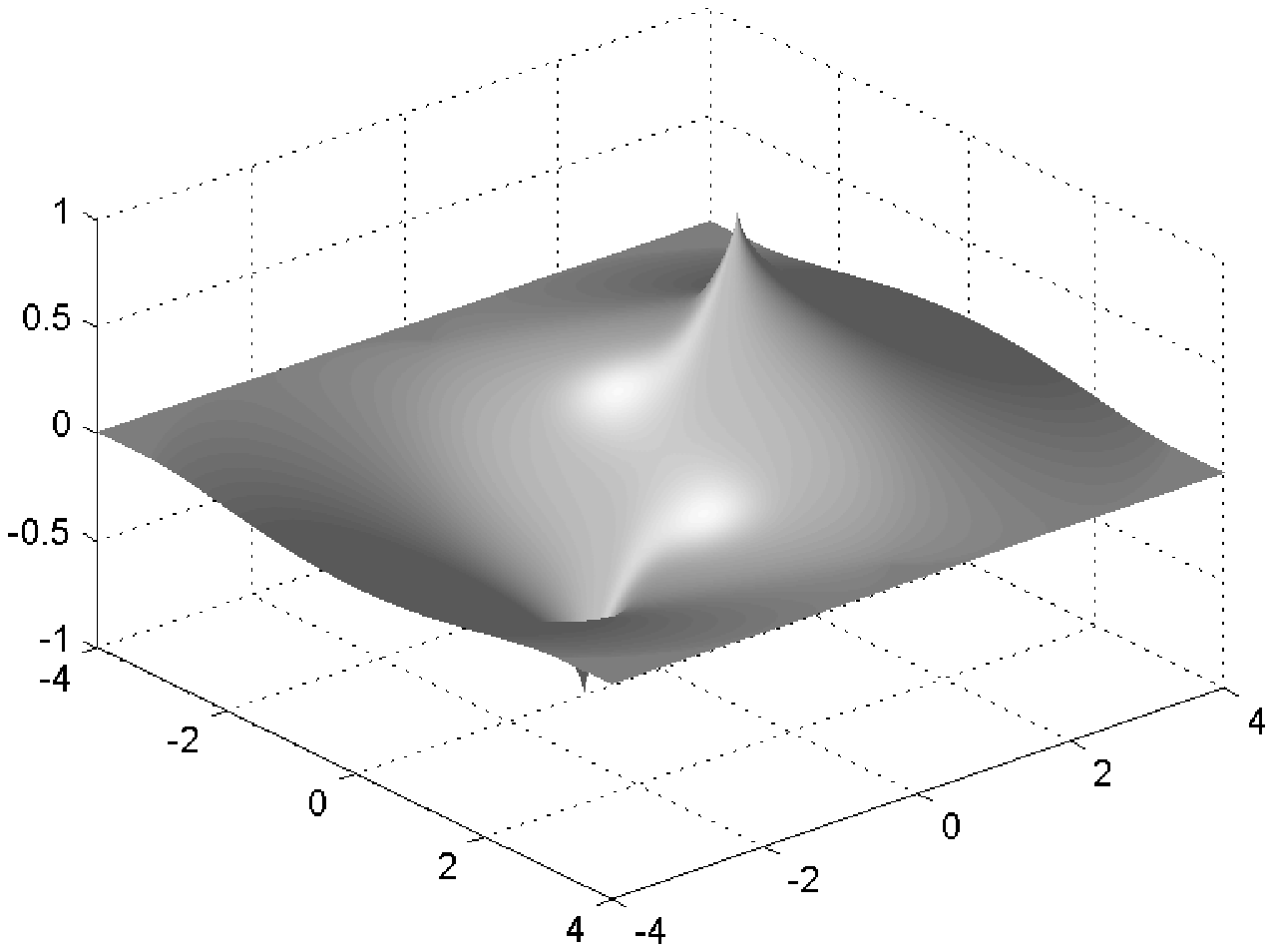}
\mbox{} \vspace{-0.850cm} \\
\mbox{} \hspace{+1.250cm}
{\small {\bf Fig.\;1:}
The solution $u$ given
in Proposition 4.1
(numerically computed). \\
\mbox{} \hspace{+1.250cm}
Here, $ n = 2 $, $ p = 4 $,
$P = \{\;\!(-1,0), \;\!(1,0)\;\!\} $,
$ u(-1,0) = -1 $, and $ u(1,0) = 1 $.
}
\end{minipage}

\mbox{} \vspace{-0.350cm} \\
\mbox{} \hspace{+0.000cm}
\begin{minipage}[t]{14.500cm}
\includegraphics[keepaspectratio=true, scale=1]{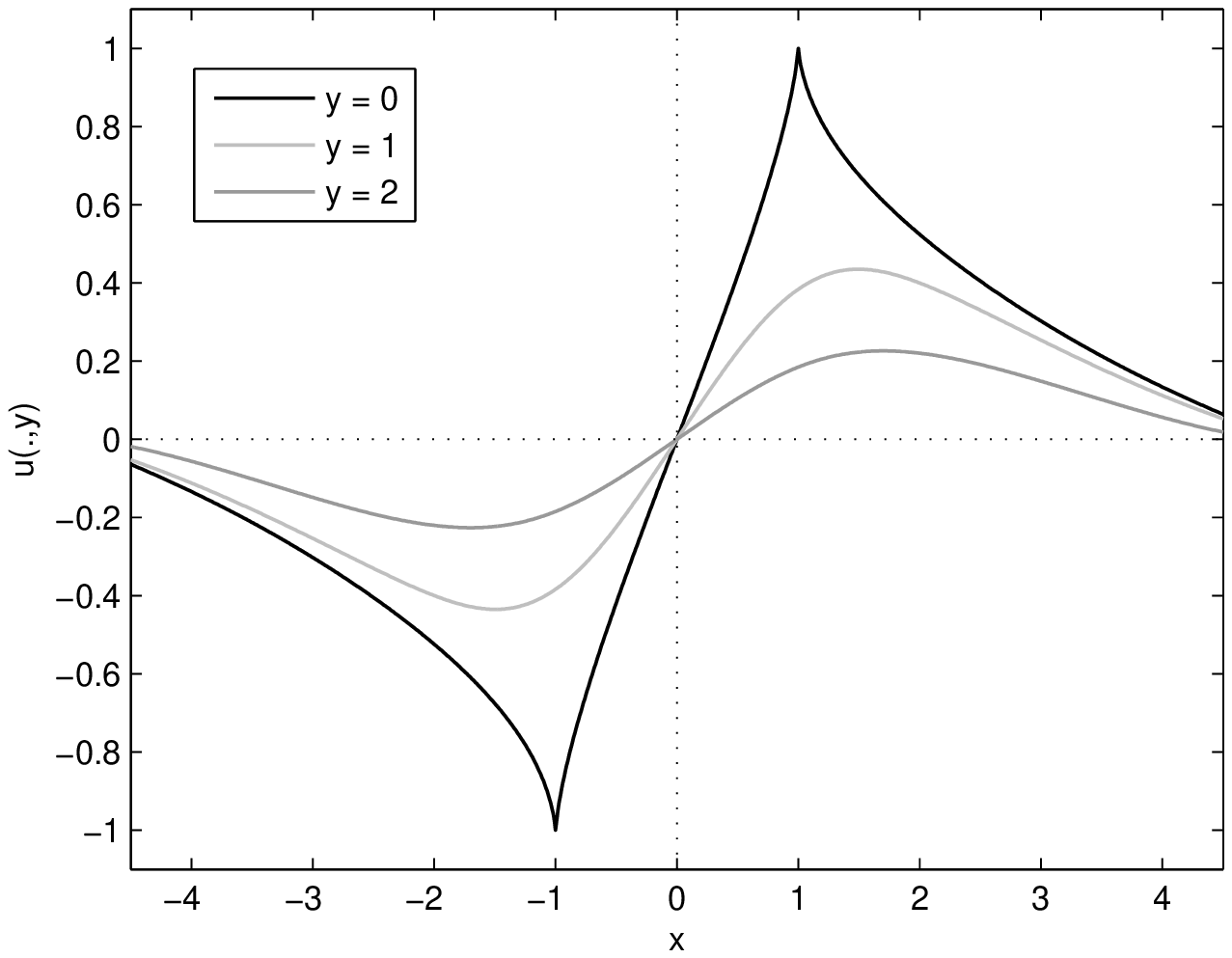}
\mbox{} \vspace{-0.500cm} \\
\mbox{} \hspace{+0.850cm}
{\small {\bf Fig.\;2:}
Three particular $x$\;\!-\;\!profiles of
the solution $u = u(x,y)$
shown above \\
\mbox{} \hspace{+0.875cm}
in Fig.\;1,
corresponding to
the sections
$ y = 0 $, $ y = 1 $,
and $ y = 2 $, respectively.
}
\end{minipage}

\mbox{} \vspace{-0.550cm} \\

\

\end{document}